\documentclass[12pt]{amsart}
\usepackage{amssymb}
\usepackage{psfig}

\theoremstyle{plain}
\newtheorem*{thm}{Theorem}
\newtheorem*{lem}{Lemma}
\newtheorem*{prop}{Proposition}
\newtheorem*{cor}{Corollary}

\theoremstyle{remark}

\newtheorem*{notat}{Notations and Conventions}

\def\hgt{\operatorname{height}}
\def\gldim{\operatorname{gldim}}
\def\GL#1{\operatorname{GL}_#1(\mathbb{Z})}
\def\co1#1#2{\operatorname{H}^1(#1,#2)}
\def\Hom{\operatorname{Hom}}
\def\End{\operatorname{End}}
\def\Id{\operatorname{Id}}

\def\Ker{\operatorname{Ker}}
\def\Cl{\operatorname{Cl}}
\def\U{\operatorname{U}}
\def\gp{\operatorname{gp}}
\def\ch{\operatorname{char}}
\def\rank{\operatorname{rank}}
\def\diag{\operatorname{diag}}
\def\GLA{\operatorname{GL}(A)}
\def\GLV{\operatorname{GL}(V)}
\def\Supp{\operatorname{Supp}}
\def\Max{\operatorname{Max}}
\def\Spec{\operatorname{Spec}}
\def\Sing{\operatorname{Sing}}
\def\Reg{\operatorname{Reg}}
\def\M{\mathfrak{M}}
\def\m{\mathfrak{m}}
\renewcommand{\P}{\mathfrak{P}}
\def\Z{\mathbb{Z}}
\def\R{\mathbb{R}}
\def\Q{\mathbb{Q}}
\def\p{\mathfrak{p}}
\def\W{\mathcal{W}}
\def\onto{\twoheadrightarrow}
\def\mono{\rightarrowtail}
\def\bij{\stackrel{1-1}{\longleftrightarrow}}

\def\iso{\stackrel{\cong}{\longrightarrow}}
\def\iff{\Longleftrightarrow}
\def\osum#1{\sigma(#1)}
\def\Bar{\overline{\phantom{x}}}
\def\bar#1{\overline{#1}}
\def\gen#1{\langle#1\rangle}

\begin{document}


\title{Multiplicative invariants and semigroup algebras}

\author{Martin Lorenz}
\address{Department of Mathematics\\
        Temple University\\
        Philadelphia, PA 19122-6094}
\email{lorenz@math.temple.edu}
\thanks{Research of the author supported in part by NSF Grant
DMS-9618521}



\begin{abstract}
Let $G$ be a finite group acting by automorphism on a lattice
$A$, and hence on the group algebra $S=k[A]$. The algebra of
$G$-invariants in $S$ is called an algebra of multiplicative
invariants.

We investigate when algebras of multiplicative
invariants are semigroup algebras. In particular, we present an
explicit version of a result of Farkas stating that multiplicative
invariants of finite reflection groups are indeed semigroup algebras.
On the other hand, multiplicative invariants arising from fixed
point free actions are shown to never be semigroup algebras.
In particular, this holds whenever $G$ has odd prime order. 
\end{abstract}

\maketitle


\section*{Introduction} \label{intro}

This article continues our investigaton of multiplicative 
invariants in \cite{L,L1,L2}
and is motivated by Farkas' work in \cite{Fa,Fa2,Fa3}. 

Our specific focus here is a suitable permanence theorem
for multiplicative actions of finite groups
analogous to the classical Shephard-Todd-Chevalley Theorem for ``linear"
actions of finite groups (of good order) on polynomial algebras;
this theorem states precisely when the corresponding algebra
of invariants is again a polynomial algebra (e.g., \cite[p.~115]{Bou}).

Multiplicative actions, also called exponential actions 
\cite{Bou}, are certain group actions on Laurent polynomial rings 
or, equivalently, group algebras of lattices.
Specifically, let $A$ denote a lattice, i.e., a free 
abelian group of finite rank, and let 
$G$ be a group acting by automorphisms on $A$. This action 
extends uniquely to an action of $G$ on the group algebra $S=k[A]$ 
of $A$. Actions of this type are
referred to as \emph{multiplicative actions}, and the resulting 
algebra of invariants $R=S^G$ is called an algebra of 
\emph{multiplicative invariants}.
It is easy to see that $R$ is again
a group algebra only if $G$ acts trivially on $A$; see 
(\ref{multinv=semigroup}). Thus the permanence theorem we have
in mind is a characterization of all multiplicative actions
yielding invariants that are semigroup algebras. Unfortunately,
this article falls short of reaching this goal.

Here is the state of affairs and our contribution.
It is implicit in \cite[proof of Theorem 10]{Fa2} that
multiplicative invariants of finite reflection groups are
indeed semigroup algebras; this has been pointed out by Farkas
himself in \cite[p.~72]{Fa3}. After deploying the requisite
background material and some technicalities in sections
\ref{semigrpalg} and \ref{multinv},
we present in section \ref{reflgps} an explicit proof of Farkas' 
theorem, along with
an analysis of the structure of the corresponding
semigroup and of the class group of the invariant algebra. 
The result, Theorem (\ref{reflinv}), is derived from 
a classical fact 
\cite[Th\'eor\`eme 1 on p.~188]{Bou} concerning multiplicative
invariants of Weyl group actions on weight lattices of root 
systems. The method employed leads directly to an explicit
fundamental system of invariants. I don't know if the
converse of Farkas' theorem holds: Do all multiplicative invariants
that are semigroup algebras come from reflection groups?
Section \ref{fpf} offers a first step in this direction.
We show in Theorem (\ref{negative}) that multiplicative
invariants of fixed point free actions (in rank at least 2)
are never semigroup algebras. In particular, this holds for
all multiplicative actions of finite groups of odd prime order.
Our main tool in this section is an investigation of the
singularities of multiplicative invariants. 
Doubtless, a good deal more can be said on this topic.
Ultimately, the complete proof of the envisioned converse
to Farkas' theorem will likely involve an investigation of certain group
actions on semigroup algebras rather than just group algebras,
and this will actually presumably be the proper setting for the
theorem.

We hope to return to these questions in a future publication.


\begin{notat}
Throughout this note, $k$ will denote a commutative field.
All monoids considered in this article are 
understood to be commutative. We use $\Z_+$ to denote the
set of nonnegative integers and similarly for $\R_+$.
Further notation will be introduced below, in particular in 
(\ref{basics}).
\end{notat}


\section{Semigroup Algebras} \label{semigrpalg}


\subsection{Commutative semigroup algebras} \label{semigroup}
Let $M$ denote a monoid, with operation written
as multiplication and identity element $1$, and let $k[M]$
denote the semigroup (or monoid) algebra of $M$ over $k$. 
Thus every element $\alpha\in k[M]$ can be uniquely written
in the form
$$
\alpha=\sum_{m\in M}k_mm\quad\text{with $k_m\in k$ almost all
zero.}
$$
The set $\{m\in M\mid k_m\neq 0\}$ is called the \emph{support}
of $\alpha$, denoted $\Supp(\alpha)$. Multiplication in $k[M]$
is defined by $k$-linear extension of the multiplication of $M$.

A good reference for general ring theoretic properties of
commutative semigroup algebras is \cite{G}.
We note in particular the following facts:
\begin{itemize}
\item The $k$-algebra $k[M]$ is finitely generated 
(\emph{affine}) if and only if
$M$ is a finitely generated monoid. This is trivial.
\item $k[M]$ is a domain iff $M$ is \emph{cancellative}
($ax=ay\Rightarrow x=y$ for $a,x,y\in M$) and
\emph{torsion-free} ($x^n=y^n\ (n>0)\Rightarrow x=y$
for $x,y\in M$); see \cite[Theorem 8.1]{G}.
\item Assume $k[M]$ is a domain. Then $k[M]$ is integrally
closed iff $M$ is
\emph{normal}: $x^n=y^nz$ for $x,y,z\in M$ implies 
$z=z_1^n$ for some $z_1\in M$; see \cite[Corollary 12.11]{G}.
\end{itemize}


\subsection{Affine normal semigroups} 
\label{affine}
Finitely generated cancellative
torsion-free normal monoids are often simply referred to as
\emph{affine normal semigroups}. 
By (\ref{semigroup}), we have for any monoid $M$:
\begin{quote}
\emph{The $k$-algebra $k[M]$ is an affine integrally closed 
domain iff $M$ is an affine normal semigroup.}
\end{quote}
As a reference for affine semigroup algebras in particular, 
I recommend \cite{BH}. By \cite[Proposition 6.1.3]{BH},  
affine normal semigroups $M$ have the following structure:
\begin{quote}
\emph{$M = \U(M)\times M_+$, where $\U(M)$ , the group of 
units of $M$, is a free abelian group of 
finite rank and $M_+$ is an affine normal
semigroup that is \emph{positive}, that is, $\U(M_+)=\{1\}$.}
\end{quote}
Further, by \cite[Theorem 11.1]{G}, the group of units of $k[M]$
is given by:
$$
\U(k[M])=k^*\times \U(M)\ .
$$
The $k$-algebra map $\mu: k[M]\to k$ that is
given by $\mu(m)=1$ for $m\in \U(M)$ and $\mu(m)=0$ if
$1\neq m\in M_+$ is called the \emph{distinguished augmentation}
of $k[M]$.


\subsection{Gubeladze's polytope} 
\label{polytope}
Let $M$ be a finitely generated monoid that is 
cancellative and torsion-free, and let $\gp(M)$ denote the
group of fractions of $M$. Our hypotheses on $M$ imply that
$\gp(M)$ is free abelian of finite rank and the natural map
$M\to\gp(M)$ is an embedding. Thus $M$ is embedded in the
real vector space $V=\R\otimes\gp(M)$. Denote by $C(M)$ the
convex cone in $V$ that is spanned by $M$; so (using additive
notation in $V$),
$$
C(M)=\R_+M=\{r_1m_1+\dots+r_km_k\mid r_i\in\R_+, m_i\in M\}\ .
$$
Assume now that $M$ is positive, that is, $\U(M)$ is trivial. 
Then there exists an affine hyperplane $H$ in $V$ with $0\notin H$
such that $C(M)$ is the convex cone in $V$ that is spanned by
$$
\Phi(M)=H\cap C(M)\ .
$$
The set $\Phi(M)$ is a polytope (the convex hull of finitely
many points) in $H$. Different choices of $H$ lead to projectively
equivalent polytopes; see \cite{Gu}, \cite{Gu2}.

Following Gubeladze \cite{Gu2}, $M$ is called 
\emph{$\Phi$-simplicial} iff $\Phi(M)$ is a simplex (the convex
hull of finitely many affinely independent points). For
several equivalent characterizations of this notion,
see \cite[Proposition 1.1]{Gu2}.


\subsection{Torus action} 
\label{torus}
Let $M$ be an affine normal semigroup. Then, as we observed
above, $M$ embeds into the lattice $L=\gp(M)\cong\Z^n$ for some $n$.
Assuming $k$ algebraically closed for now, the algebraic torus
$T=\Hom(L,k^*)\cong \left(k^*\right)^n$ operates on $k[M]$
via
$$
m^{\tau}=\tau(m)^{-1}m\qquad(\tau\in T, m\in M)
$$
and $k$-linear extension to all of $k[M]$. This operation is
rational, and the corresponding operation on 
$\Max k[M]=\Hom(k[M],k)$ has the distinguished augmentation
as its only fixed point if $\U(M)$ is trivial, and no
fixed points at all otherwise. For further background
on the geometric aspect of affine normal semigroup algebras,
see \cite{Fu}; the above assertion about fixed points is
an (easy) exercise in \cite[p.~28]{Fu}.


\section{Multiplicative Invariants}\label{multinv}


\subsection{Basics} 
\label{basics}
The following notation will be kept throughout this article:
\bigskip
\begin{tabbing}
\hspace*{.35in}\=\hspace{1in}\=\kill
\> $A$ \> will be a free abelian group of finite rank;\\
\> $S=k[A]$ \> will denote the group algebra of $A$ over $k$;\\
\> $G$ \> will be a finite group acting be automorphisms on $A$,\\
\> \> and hence on $S$ as well; the action will be written\\
\> \> exponentially, $a\mapsto a^g$;\\
\> $R=S^G$ \> is the subalgebra of $G$-invariants in $S$.
\end{tabbing}
\bigskip
In this situation, $A$ is often called a \emph{$G$-lattice}. 
As our main concern is $R$, 
the algebra of multiplicative $G$-invariants, we may
assume that the $G$-lattice $A$ is \emph{faithful}, 
that is, the map $G\to\GLA$ that defines the $G$-action is 
injective. Finally, $A$ will be called \emph{effective}
if the subgroup $A^G$ of $G$-invariant elements of $A$
is trivial.

The \emph{orbit sum} of an element $a\in A$ is the element 
of $S$ that is defined by
$$
\osum{a}=\sum_{x\in a^G}x\ \in S\ .
$$
where $a^G=\{a^g\mid g\in G\}\subseteq A$ denotes the $G$-orbit 
of $a$. Orbit sums are clearly $G$-invariant, and hence they actually 
belong to $R$. In fact, they provide a $k$-basis
for $R$ :
$$
R=\bigoplus_{a\in A/G}k\osum{a}\ ,
$$
where $A/G$ denotes a transversal for the $G$-orbits in $A$.
As $k$-algebra, $R$ is an affine integrally closed domain;
all these properties are inherited from $S$. (Note that $A$ 
is an affine normal semigroup.)


\subsection{Passage to an effective lattice} 
\label{effective}
Let $\Bar$ denote the canonical map $A\onto A/A^G$ and its 
extension to $S$; so
$$
\Bar : S=k[A] \onto \bar{S}=k[A/A^G]\ ,\quad a\mapsto aA^G\ 
(a\in A) \ .
$$
Note that $\bar{A}=A/A^G$ is a $G$-lattice and the map
$\Bar$ is $G$-equivariant. Moreover, letting $G_x$ denote
the isotropy (stabilizer) subgroup of $G$ of an element
$x$ in $A$ or in $\bar{A}$, we have
\begin{equation}\label{E:isotropy}
G_a=G_{\bar{a}}\quad\text{for all $a\in A$.}
\end{equation}
Here, the inclusion $G_a\subseteq G_{\bar{a}}$ is clear. The
reverse inclusion follows from the fact that the map
$G_{\bar{a}}\to A^G$, $g\mapsto a^ga^{-1}$, is a group
homomorphism, and hence it must be trivial, as $G_{\bar{a}}$
is finite while $A^G$ is torsion free.
We deduce from the above equality of isotropy groups that
\begin{quote}
\emph{$\bar{A}$ is an effective $G$-lattice.}
\end{quote} 
Further, $\Bar:S\to\bar{S}$ sends the orbit sum
$\osum{a}$ to the orbit sum $\osum{\bar{a}}$, and
$\osum{\bar{a}}=\osum{\bar{b}}$ is equivalent with
$\osum{a}=\osum{b}c$ for some $c\in A^G$.
Consequently, 
\begin{quote}
\emph{The map $\Bar$ maps $R$ \emph{onto} the $G$-invariants
in $\bar{S}$, that is, $\bar{R}=\bar{S}^G$. The kernel of this
epimorphism is the ideal $\left(a-1\mid a\in A^G\right)$ of $R$.}
\end{quote} 
Finally, every $G$-equivariant homomorphism from $A$ to 
some effective $G$-lattice clearly factors through 
$\Bar:A\to\bar{A}$.


\subsection{Multiplicative invariants that are semigroup 
algebras} 
\label{multinv=semigroup}
In this section, we note some consequences of the
assumption that $R$ is a semigroup algebra. In particular,
it will turn out that $\bar{R}$ is a semigroup algebra as
well in this case.

\begin{prop} Assume that $\varphi:k[M]\iso R$ for some 
semigroup algebra $k[M]$.
Then $M$ is an affine normal semigroup, so $M=\U(M)\times M_+$
as in \textnormal{(\ref{affine})}, and $M_+$ is 
$\Phi$-simplicial.
Moreover, the isomorphism $\varphi$ can be chosen so that 
$\varphi(\U(M))=A^G$.  Finally, $\varphi$ restricts to an 
isomorphism $k[M_+]\iso\bar{R}$, in the notation of 
\textnormal{(\ref{effective})}.
\end{prop}

\begin{proof} 
First, $M$ must be an affine normal semigroup, 
since $R$ is an affine integrally closed domain; see
(\ref{semigroup}) and (\ref{basics}). 
Hence, $M = \U(M)\times M_+$ and
$$
k^*\times \U(M)=\U(k[M])\cong\U(R)=\U(S)^G=
k^*\times A^G\ .
$$
Letting $\alpha: S\to k$ denote the
distinguished augmentation of $S=k[A]$, sending all $a\in A$
to $1$, the given isomorphism $\varphi$
can be modified by defining 
$\psi(m)=\alpha\varphi(m)^{-1}\varphi(m)$
for $m\in\U(M)$ and $\psi(m)=\varphi(m)$ for $m\in M^+$ to 
obtain
a new isomorphism $\psi: k[M]\iso R$ which maps $\U(M)$ 
onto $A^G$. 
The composite $\Bar\circ\psi: k[M]\iso R\onto\bar{R}$ has 
kernel 
$\psi^{-1}\left((a-1\mid a\in A^G)\right)=(m-1\mid m\in\U(M))$. 
Hence, this map restricts to an isomorphism $k[M_+]\iso\bar{R}$.

It remains to show that $M_+$ is $\Phi$-simplicial. 
Now, $M_+$ is $\Phi$-simplicial if and only if
$k[M_+]$ is almost factorial, that is, the class group
$\Cl(k[M_+])$ is torsion; see \cite[Proposition 1.6]{Gu2}.
However, $\Cl(k[M_+])\cong\Cl(k[M])$, by
\cite[Corollary 7.3 and Theorem 8.1]{F}. Therefore, 
$\Cl(k[M])\cong\Cl(R)$. Since $\Cl(S)$ is trivial, it follows
form Samuel's
theory of Galois descent (cf.~\cite[Theorem 16.1]{F})
that $\Cl(R)$ embeds into 
$\co1{G}{\U(S)}=\Hom(G,k^*)\oplus\co1{G}{A}$, a finite 
$|G|$-torsion group.
(The precise form of class groups of multiplicative
invariants is known \cite{L}, but this information is not
needed here.) We deduce that $\Cl(k[M_+])$ is finite
$|G|$-torsion, thereby completing the proof.
\end{proof}

As a very simple special case, assume that $R$ is actually 
a group algebra; so $M=\U(M)$. Then the above proposition 
yields that $R=k[A^G]$. Now $S$ is integral over $R=S^G$ and, 
on the other hand, $A/A^G$ is torsion-free. Thus we must 
have $A=A^G$, whence $G$ acts trivially on $A$. This 
substantiates a remark made in the introduction.


\subsection{A reduction lemma} 
\label{reduction}
In this section, we will prove a technical lemma stating
that an algebra of multiplicative invariants is a
semigroup algebra provided a closely related one is. Let
$$
\mathcal{M}(A)
$$ 
denote the submonoid of $(R,\cdot)$ that is
generated by the orbit sums $\osum{a}$ for 
$a\in A$, and similarly for other $G$-lattices. 

\begin{lem} Let $A\subseteq B$ be $G$-lattices such that
$B/A$ is $G$-trivial. Suppose that $k[B]^G=kC$, the $k$-linear
span of some subset $C\subseteq\mathcal{M}(B)$. Then $k[A]^G=kD$ with
$D=C\cap k[A]$.
\end{lem}

\begin{proof}
Note that $D$ is a subset of $k[A]^G$; so clearly $kD\subseteq
k[A]^G$. For the other inclusion, let $\alpha\in k[A]^G$ be
given. Then $\alpha=\sum_{c\in C}k_cc$ with $k_c\in k$ almost all
zero. We show by induction on the minimum number, $n(\alpha)$,
of nonzero terms in such an expression that $\alpha\in kD$.
The case $n(\alpha)=0$ (i.e., $\alpha=0$) being obvious, assume
$\alpha\neq 0$. Then some $d\in C$ with $k_d\neq 0$ must
satisfy $\Supp(d)\cap A\neq\emptyset$. Say 
$d=\osum{b_1}\cdot\ldots\cdot\osum{b_l}$ with $b_j\in B$. Then
$$
\Supp(d)\subseteq\{b_1^{g_1}\cdot\ldots\cdot b_l^{g_l}
\mid g_j\in G\}\ .
$$
So some product $b_1^{g_1}\cdot\ldots\cdot b_l^{g_l}$ belongs to
$A$. Inasmuch as $B/A$ is $G$-trivial, all these products are
congruent to each other modulo $A$, and hence they all belong 
to $A$. Thus, $\Supp(d)\subseteq A$ and so $d\in D$. Since
$\alpha-k_dd$ belongs to $kD$, by induction, we conclude that
$\alpha\in kD$ as well. This proves the lemma.
\end{proof}

Note that if the subset $C$ in the Lemma is $k$-independent or 
multiplicatively closed then so is $D=C\cap k[A]$.
Hence, if $k[B]^G=kC$ is a semigroup algebra, with semigroup
basis $C$, then $k[A]^G=kD$ is a semigroup algebra with semigroup 
basis $D$. 

We also remark, for future use, that  
the argument in the proof of the Lemma shows that, for 
$d=\prod_{j=1}^l\osum{b_j}\in\M(B)$,
\begin{equation}\label{E:product}
\prod_{j=1}^l\osum{b_j}\in k[A]\iff\Supp(d)\cap a\neq\emptyset
\iff\prod_{j=1}^l b_j\in A\ .
\end{equation}


\section{Reflection Groups}\label{reflgps}


\subsection{Reflections}\label{refl}
An endomorphism $\phi$ of a vector space is called a 
\emph{pseudoreflection} if $\Id-\phi$ has rank 1; 
$\phi$ is a \emph{reflection} if, in addition, $\phi^2=\Id$.
 
Keeping the notation of (\ref{basics}), we will assume in this
section that $A$ is a $G$-lattice which, without 
essential loss, will be assumed faithful. We will 
further assume that $G$ is a reflection group on $A$; so:
\begin{quote}
\emph{$G$ is a finite subgroup of $\GLA$ that is generated
by reflections.}
\end{quote}
Here, an element $g\in G$ is called a \emph{reflection} if $g$
is a reflection on $A\otimes_{\Z}\Q$. We remark that, since
$\det g=\pm 1$ holds for all $g\in G$, pseudoreflections in
$G$ are automatically reflections. 
They can also be characterized by the condition that the
subgroup $A^{\gen{g}}=\Ker_A(g-\Id)$ of $g$-fixed points 
in $A$ have rank equal to $\rank(A)-1$ or, alternatively,
\begin{quote}
\emph{$\Ker_A(g+\Id)=\{a\in A\mid a^g=a^{-1}\}$ is infinite cyclic.}
\end{quote}
As in (\ref{effective}), we let $\Bar$ denote the canonical map 
$A\onto\bar{A}=A/A^G$. Note that (\ref{effective})(\ref{E:isotropy})
implies that $\bar{A^{\gen{g}}}=\bar{A}^{\gen{g}}$ holds for all 
$g\in G$.
Therefore, if $g$ acts as a reflection on $A$ then it does so on 
$\bar{A}$ as well, and conversely. 


\subsection{Root systems}\label{roots}
Embed $A$ into the $\R$-vector space $V=A\otimes_{\Z}\R$ 
and view $G$ as a subgroup of $\GLV$.
As is customary, we 
will use additive notation in $A$ and $V$.
Define
$$
\rho(v)=|G|^{-1}\sum_{g\in G}v^g\qquad (v\in V)\ .
$$
Thus, $\rho$ is an idempotent $\R[G]$-endomorphism of $V$
with $\rho(V)=V^G$, the subspace of $G$-fixed points in $V$.
Putting $\pi=1-\rho\in\End_{\R[G]}(V)$, we obtain
$$
A\subseteq\rho(A)\oplus\pi(A)\subseteq\rho(V)\oplus\pi(V)
=V\ .
$$

For each reflection $g\in G$, let the two possible
generators of $\Ker_A(g+\Id)$ be denoted $\pm a_g$. Define
$$
\Phi=\Phi_{A,G}=\{\pm a_g\mid 
\text{$g$ a reflection in $G$}\}\ .
$$
The crucial properties of $\Phi$ are listed in the following
lemma due to Farkas \cite[Lemmas 1--3]{Fa2}.
\begin{lem}
$\Phi=\Phi_{A,G}$
is a reduced crystallographic root system in $\pi(V)$, and 
the restriction of $G$ to $\pi(V)$ is the Weyl group $\W(\Phi)$
of $\Phi$. Furthermore, 
$$
\Z\Phi\subseteq A\subseteq\pi^{-1}(\Lambda)\ ,
$$ 
where $\Z\Phi$, the $\Z$-span of $\Phi$ in $V$, is the 
\emph{root lattice} and $\Lambda=\Lambda_{A,G}=\{v\in\pi(V)\mid
v-v^g\in{\Z}a_g\text{ for all reflections $g\in G$}\}$ is
the \emph{weight lattice} of $\Phi$.
\end{lem}

For background on root systems, we refer to \cite{Bou} or \cite{H}.


\subsection{Multiplicative invariants of reflection groups}
\label{reflinv}
Our goal here is to prove the following result implicit in the
work of Farkas \cite{Fa2,Fa3}. We will use the notation of
(\ref{roots}). 

\begin{thm} Let $A$ be a free abelian group of finite rank,
and let $G$ be a finite subgroup of $\GLA$ that is generated
by reflections. Then the invariant algebra $R=k[A]^G$ is a
semigroup algebra; in fact, $R\cong k[M]$ with 
$M=A^G\times (\pi(A)\cap\Lambda_+)$, where $\Lambda_+$
is the semigroup of dominant weights for some
base of the root system $\Phi_{A,G}$.
\end{thm}

\begin{proof} Fix a base $\Delta=\{\alpha_1,\ldots,\alpha_r\}$ 
for $\Phi=\Phi_{A,G}$, i.e., $\Delta$ is a subset of $\Phi$
that is an $\R$-basis of $\pi(V)$ and such that 
$\Phi\subseteq\Z_+\Delta\cap-\Z_+\Delta$.
So $\alpha_i=\pm a_{g_i}$ for certain reflections $g_i\in G$.
The \emph{fundamental dominant weights} 
$\lambda_1,\ldots,\lambda_r$ are determined by 
$\lambda_i-\lambda_i^{g_j}=\delta_{i,j}\alpha_j$ (Kronecker delta); 
they form a $\Z$-basis of the weight lattice $\Lambda$. The
semigroup $\Lambda_+$ of dominant weights for $\Delta$ is
$$
\Lambda_+=\oplus_{i=1}^r\Z_+\lambda_i\ .
$$
It is a classical result \cite[Th\'eor\`eme 1 on p.~188]{Bou}
that $k[\Lambda]^G$ is a polynomial algebra, with the orbit sums
of the fundamental dominant weights as independent generators.
In other words, 
\begin{quote}
\emph{$k[\Lambda]^G=kE$, with 
$E=\gen{\osum{\lambda_1},\ldots,\osum{\lambda_r}}\cong\Lambda_+$
a $k$-independent submonoid of $\mathcal{M}(\Lambda)$.}
\end{quote}
Now put $B=\rho(A)\oplus\Lambda$, a $G$-lattice in $V$ with
$A\subseteq B$ and $B/A$ $G$-trivial.
To see the latter, note that $A$ contains $A^G\oplus\Z\Phi$,
and $B/(A^G\oplus\Z\Phi)\cong 
\left(\rho(A)/A^G\right)\oplus\left(\Lambda/\Z\Phi\right)$ is
$G$-trivial, since $\rho(A)\subseteq V^G$ and the Weyl group
$G$ of $\Phi$ acts trivially on the fundamental group 
$\Lambda/\Z\Phi$ of $\Phi$; cf.~\cite[p.~167]{Bou}.
Inasmuch as $k[B]=k[\rho(A)]\otimes_kk[\Lambda]$, 
with $\rho(A)=B^G$, the $G$-invariants in $k[B]$ are  
given by $k[B]^G=k[B^G]\otimes_kk[\Lambda]^G$.
Thus, using the above description of $k[\Lambda]^G$,
$$
k[B]^G = k[B^G]\otimes_kkE = kC\quad\text{\emph{with }}
\ C=B^G\times E\ .
$$
Note that $C$ is a $k$-independent submonoid of $\mathcal{M}(B)$. 
Lemma (\ref{reduction}) therefore implies that
$k[A]^G=kD$ is a semigroup algebra, with semigroup basis
$D=C\cap k[A]$. It remains to verify the description
of the monoid given in the theorem.
To this end, note that, by (\ref{reduction})(\ref{E:product}), 
the
isomorphism $B^G\oplus \Lambda_+\iso B^G\times E=C$ restricts 
to an isomorphism $M:=(B^G\oplus\Lambda_+)\cap A\iso D$. 
Furthermore, writing $a\in A$ as $a=\rho(a)+\pi(a)$, we see
that $a\in M$ if and only if 
$\pi(a)\in\Lambda_+$. Since $\Ker_A(\pi)=A^G$ and $\bar{A}=A/A^G$
is free, we have $A=A^G\oplus A'$ with $A'\cong \pi(A)$ via $\pi$.
This decomposition induces a corresponding one for $M$, because
$A^G\subseteq M$; so $M=A^G\oplus (M\cap A')$ and $M\cap A'\cong
\pi(A)\cap\Lambda_+$ via $\pi$.
This completes the proof of the theorem.
\end{proof}


\subsection{Generators} 
\label{generators}
We now descibe how the foregoing leads to an explicit set of
\emph{fundamental invariants}, that is, algebra generators
for $R$. Inasmuch as $R\cong k[M]$, this amounts to finding
generators for $M$ and tracing them through the isomorphism.
As this isomorphism is the identity on $\U(M)=A^G$,
we will concentrate on $M_+$.

\subsubsection{Generators for $M_+=\pi(A)\cap\Lambda_+$}
\label{Mgenerators}
Since the semigroup $M_+$ is positive, it has a unique
minimal generating set, the so-called \emph{Hilbert basis} of
$M_+$.
Here, in outline, is how to find this Hilbert basis; 
for complete
details and an algorithmic treatment, see \cite[Chapter 13]{St}.

Recall that $\Lambda_+=\oplus_{i=1}^r\Z_+\lambda_i$, where
$\lambda_1,\ldots,\lambda_r$ are the fundamental dominant weights.
These belong to $\pi(A)\otimes\Q\subseteq V$. Hence, there are
suitable $0\neq z_i\in\Z_+$ so that $m_i=z_i\lambda_i\in M_+$; we
will assume that $z_i$ is chosen minimal.
The subset
$$
K=\sum_{i=1}^r[0,m_i]=
\{\sum_{i=1}^r t_im_i\mid 0\le t_i\le 1\}
$$
of $V$ is compact (a zonotope), and hence its intersection 
$K\cap M_+$ with
the discrete $M_+$ is finite. It is easy to see that $K\cap M_+$
generates $M_+$; the Hilbert basis of $M_+$ can be found by
selecting the indecomposable elements of $K\cap M_+$, that is,
the elements $m\in K\cap M_+$ that cannot be written as 
$m=n+n'$ with $0\neq n,n'\in K\cap M_+$. Note that $m_1,\ldots,m_r$
are certainly indecomposable, by the minimal choice of the $z_i$'s
and linear independence of the $\lambda_i$'s. The remaining 
indecomposables in $K\cap M_+$ (if any) will be denoted 
$m_{r+1},\ldots,m_s$; so $s\ge r=\rank(\bar{A})$.

We remark in passing that Gubeladze's polytope $\Phi(M_+)$ is
the convex hull of $m_1,\ldots,m_r$ (up to projective equivalence;
see (\ref{polytope})). Indeed, since $m_i\in M_+=\Z_+(K\cap M_+)
\subseteq \R_+\{m_1,\ldots,m_r\}$, we have 
$C(M_+)=\R_+\{m_1,\ldots,m_r\}$. 

\subsubsection{Fundamental invariants}\label{fundamental}
As all $m_i\in\Lambda_+=\oplus_{j=1}^r\Z_+\lambda_j$, they have a
unique representation of the form $m_i=\sum_j z_{i,j}\lambda_j$
with $z_{i,j}\in\Z_+$. For $i\le r$, this representation is
simply $m_i=z_i\lambda_i$, as above. Thus we obtain the following
system of fundamental invariants:
$$
\mu_i=\prod_{j=1}^r\osum{\lambda_i}^{z_{i,j}}\ (i=1,\dots,s)
$$
Here, $\mu_1=\osum{\lambda_1}^{z_1},\ldots,
\mu_r=\osum{\lambda_r}^{z_r}$ are algebraically independent,
as the $\osum{\lambda_i}$'s are, and $R$ is a finite module
over the polynomial algebra $k[\mu_1,\ldots,\mu_r]$, since each
$\mu_i$, raised to a suitable power, belongs to 
$\gen{\mu_1,\ldots,\mu_r}$. In fact, since $R$ is Cohen-Macaulay
(cf.~\cite[Theorem 6.3.5]{BH}),
$R$ is a free module over $k[\mu_1,\ldots,\mu_r]$.


\subsection{The class group}\label{classgroup}
The formula given in \cite{L} for the class group of $R$ can be
rewritten in terms of the above root system data. Indeed, by 
(\ref{multinv=semigroup}), $R=k[M]=k[\U(M)]\otimes\bar{R}$ is
a Laurent polynomial extension of $\bar{R}$, and so $\Cl(R)=\Cl(\bar{R})$.
Further, by \cite{L}, $\Cl(\bar{R})=\co1{G}{\bar{A}^D}$,
where $D$ denotes the subgroup of $G$ that is generated by those 
reflections that are \emph{diagonalizable} on $\bar{A}$, that is, 
with respect to a suitable $\Z$-basis of $\bar{A}$, they have the
form $\diag(-1,1,\ldots,1)$. Now $G$ acts as a reflection group on
$\bar{A}^D$, and the $G$-lattice $\bar{A}^D$ is effective,
as $\bar{A}$ is. Thus, \cite[Proposition 2.2.25]{Lem} gives 
$\co1{G}{\bar{A}^D}\cong\Lambda_{\bar{A}^D,G}/\bar{A}^D$. Hence,
$$
\Cl(R)\cong\Lambda_{\bar{A}^D,G}/\bar{A}^D\ .
$$
It is perhaps worth noting that 
$\Lambda_{\bar{A}^D,G}/\bar{A}^D$ is always a direct 
summand
of $\Lambda_{A,G}/\pi(A)=\Lambda_{\bar{A},G}/\bar{A}$. This follows
from the fact that $\bar{A}^D$ is a direct summand of $\bar{A}$ as 
$G$-lattices; see \cite[Lemma 2.4]{L}. 

In the special case where $A$ is effective at the outset and $G$
contains no diagonalizable reflections, the above formula
simplifies to
$$
\Cl(R)\cong\Lambda/A\ ,
$$
with $\Lambda=\Lambda_{A,G}$ as before. 

Finally, we remark that the Picard group of $R$ is trivial, as is 
in fact 
the full projective class group $K_0(R)/\gen{[R]}$. This is a 
consequence of Gubeladze's theorem \cite{Gu} stating that all 
projective modules over $R=k[M]$ are free.


\subsection{Examples}\label{examples}
We illustrate the foregoing with a couple of explicit examples.
In each case, $A$ will be effective; so $\pi=\Id$ and 
$M=M_+=A\cap\Lambda_+$. We will follow the notations in the
proof of Theorem (\ref{reflinv}) and in (\ref{generators})
quite closely.

\subsubsection{An example in rank 2} Let $A$ be free abelian of
rank 2, with $\Z$-basis $\{a,b\}$, and let $G$ be the subgroup of
$\GLA=\GL2$ that is generated by the matrices
$r=\left(
\begin{smallmatrix}
0 & 1 \\
1 & 0
\end{smallmatrix}
\right)$
and 
$s=\left(
\begin{smallmatrix}
1 & -1 \\
0 & -1
\end{smallmatrix}
\right)$. 
(These matrices act on the right on $A$, viewed as
integer row vectors of length 2.) The generators $r$ and $s$
are reflections, and $G\cong S_3$, the symmetric group on
3 symbols. The only other reflection in $G$ is 
$t=\left(
\begin{smallmatrix}
-1 & 0 \\
-1 & 1
\end{smallmatrix}
\right)$;
all reflections are conjugate in $G$, and none are diagonalizable.
As a generator for $\Ker_A(g+\Id)$, we choose 
$a_r=(-1,1)=a^{-1}b$;
similarly, we select $a_s=(0,1)=b$ for $s$ and $a_t=(1,0)=a$
for $t$. So $\Phi=\{\pm a_r, \pm a_s, \pm a_t\}$ (a root system
of type $A_2$). As base for $\Phi$, we fix 
$\Delta=\{\alpha_1=-a_t=(-1,0),\alpha_2=a_s=(0,1)\}$;
so $g_1=t$ and $g_2=s$.
This leads to the fundamental dominant weights
$\lambda_1=(-2/3,1/3)$, $\lambda_2=(-1/3,2/3)$. The
zonotope $K=[0,m_1]+[0,m_2]$ of (\ref{Mgenerators}) is 
given by $m_i=3\lambda_i$, 
and we obtain the following
generators for $M$: $m_1$, $m_2$, and $m_3=\lambda_1+\lambda_2$.
\bigskip\newline

\centerline{\psfig{figure=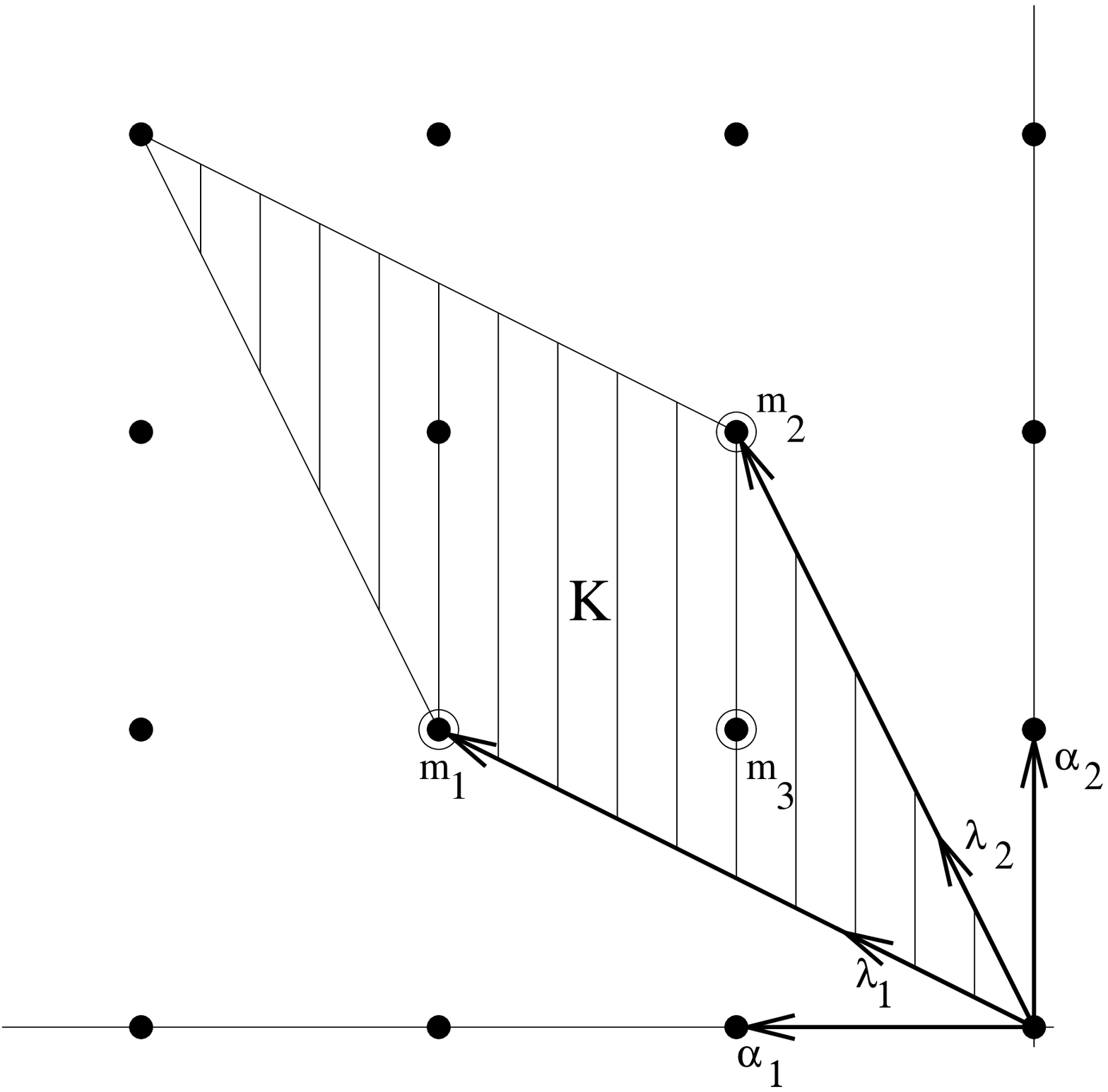,height=3in}}
\bigskip

Therefore, $\osum{\lambda_1}^3$, $\osum{\lambda_2}^3$, and
$\osum{\lambda_1}\osum{\lambda_2}$, form a fundamental 
system of invariants in $R$.
Returning to multiplicative notation, the
orbit sums for the fundamental dominant weights are:
\begin{align*}
\osum{\lambda_1} &= a^{-2/3}b^{1/3}+a^{1/3}b^{-2/3}+a^{1/3}b^{1/3}
= a^{1/3}b^{1/3}(a^{-1}+b^{-1}+1) \\
\osum{\lambda_2} &= a^{-1/3}b^{-1/3}+a^{-1/3}b^{2/3}+a^{2/3}b^{-1/3}
= a^{-1/3}b^{-1/3}(a+b+1) \ .
\end{align*}
This leads to the following explicit system of fundamental invariants:
\begin{align*}
\mu_1=\osum{\lambda_1}^3 &= ab(a^{-1}+b^{-1}+1)^3, \\
\mu_2=\osum{\lambda_2}^3 &= a^{-1}b^{-1}(a+b+1)^3,\\
\mu_3=\osum{\lambda_1}\osum{\lambda_2}&= (a+b+1)(a^{-1}+b^{-1}+1)\ . 
\end{align*}
The class group of $R$ evaluates to $\Cl(R)=\Lambda/A\cong\Z/3\Z$.


\subsubsection{Example in rank 3} Let $A$ be free abelian 
with $\Z$-basis $\{a,b,c\}$, and let $G$ be the subgroup of
$\GLA=\GL3$ that is generated by the matrices
$r=\left(
\begin{smallmatrix}
0 & 1 & 0\\
1 & 0 & 0\\
0 & 0 & 1
\end{smallmatrix}
\right)$,
$s=\left(
\begin{smallmatrix}
1 & 0 & -1 \\
0 & 1 & -1 \\
0 & 0 & -1
\end{smallmatrix}
\right)$, and
$t=\left(
\begin{smallmatrix}
0 & 0 & 1\\
0 & 1 & 0\\
1 & 0 & 0
\end{smallmatrix}
\right)$.
This group is isomorphic to $S_4$. The generators are reflections;
they are all conjugate. The complete set reflections is the
full $G$-conjugacy class: $\{r,s,t,w=r^t,u=s^t,v=s^w\}$; none are
diagonalizable.
The root system $\Phi=\Phi_{A,G}$ evaluates to
$\Phi=\{\pm(1,0,0),\pm(1,0,-1),\pm(1,-1,0),\pm(0,1,0),\pm(0,0,1),
\pm(0,1,-1)\}$. A suitable base of $\Phi$ is 
$\Delta=\{\alpha_1=a_t=(-1,0,1),\alpha_2=a_r=(1,-1,0),
\alpha_3=a_s=(0,0,-1)\}$; so $g_1=t$, $g_2=r$, $g_3=s$.
This results in the following fundamental dominant weights:
$\lambda_1=(-1/2,-1/2,1/2)$, $\lambda_2=(1/4,-3/4,1/4)$, and
$\lambda_3=(-1/4,-1/4,-1/4)$. The zonotope $K$ is spanned
by $m_1=2\lambda_1$, $m_2=4\lambda_2$, and $m_3=4\lambda_3$,
and the generators of M are: $m_1$, $m_2$, $m_3$, 
$m_4=\lambda_2+\lambda_3$, $m_5=\lambda_1+2\lambda_2$, and
$m_6=\lambda_1+2\lambda_3$.
Calculating the orbits sums:
\begin{align*}
\osum{\lambda_1} &= a^{-1/2}b^{-1/2}c^{-1/2}(a+b+c+ab+ac+bc),\\
\osum{\lambda_2} &= a^{1/4}b^{1/4}c^{1/4}(a^{-1}+b^{-1}+c^{-1}+1),\\
\osum{\lambda_3} &= a^{-1/4}b^{-1/4}c^{-1/4}(a+b+c+1).
\end{align*}
This leads to the following explicit system of fundamental invariants:
\begin{align*}
\mu_1=\osum{\lambda_1}^2 &= a^{-1}b^{-1}c^{-1}(a+b+c+ab+ac+bc)^2, \\
\mu_2=\osum{\lambda_2}^4 &= abc(a^{-1}+b^{-1}+c^{-1}+1)^4,\\
\mu_3=\osum{\lambda_3}^4 &= a^{-1}b^{-1}c^{-1}(a+b+c+1)^4,\\
\mu_4=\osum{\lambda_2}\osum{\lambda_3}&= (a+b+c+1)(a^{-1}+b^{-1}+c^{-1}+1),\\
\mu_5=\osum{\lambda_1}\osum{\lambda_2}^2 &= (a+b+c+ab+ac+bc)(a^{-1}+b^{-1}+c^{-1}+1)^2,\\
\mu_6=\osum{\lambda_1}\osum{\lambda_3}^2 &= 
(a^{-1}+b^{-1}+c^{-1}+a^{-1}b^{-1}+a^{-1}c^{-1}+b^{-1}c^{-1})(a+b+c+1)^2.
\end{align*}
For the class group of $R$, we obtain $\Cl(R)=\Lambda/A\cong\Z/4\Z$.

The calculations for this example were performed with 
{\sf GAP} (version 3.4) \cite{Sch}; 
the code is available under 
\texttt{http://www.math.temple.edu/\~{}lorenz/semigroup.html}.


\section{Fixed-Point-Free Actions}\label{fpf}

\emph{We continue with the notation of (\ref{basics}). In addition, we
assume in this section that $\ch k$ does not divide the order of $G$.
Finally, we will continue to assume that $G$ acts faithfully on $A$;
so $G\subseteq\GLA$.}


\subsection{Cotangent spaces}\label{cotangent}
The (Zariski) \emph{cotangent space} at the maximal ideal
$\M$ of $S$ is the $k$-space
$$
\M/\M^2\ .
$$
It is a $k[G^T(\M)]$-module, where 
$G^T(\M)=\{g\in G\mid s-s^g\in\M\text{ for all $s\in S$}\}$
is the inertia group of $\M$. 

\begin{lem} Assume $k$ algebraically closed. Then, for each 
maximal ideal $\M$ of $S$, there is an isomorphism of
$k[G^T(\M)]$-modules $A\otimes_{\Z}k\iso \M/\M^2$.
An element $g\in G^T(\M)$ is a pseudoreflection on $\M/\M^2$
if and only if $g$ is a reflection on $A$.
\end{lem}

\begin{proof}
Let $\mu$ be the $k$-algebra homomorphism $\mu: S \onto S/\M\iso k$.
The $k[G^T(\M)]$-isomorphism is given by
\begin{alignat*}{3}
A&\otimes_{\Z}k\ &&\iso&\M&/\M^2\\
a&\otimes 1 &&\longmapsto&\ \ \mu(a)^{-1}a&-1+\M^2\qquad(a\in A)\ .
\end{alignat*}
Clearly, an element $g\in G^T(\M)$ is a pseudoreflection on $\M/\M^2$
if and only if $g$ acts as a pseudoreflection on $A\otimes_{\Z}k$.
Thus, it suffices to show that $g\in G$ acts as a (pseudo)reflection
on $A$ if and only if $g$ does so on $V=A\otimes_{\Z}k$.
This is a consequence of the following more general equality for 
$g$-fixed point sets: 
\begin{equation}\label{E:fixed}
V^{\gen{g}}=A^{\gen{g}}\otimes_{\Z}k \text{ holds for all $g\in G$.}
\end{equation}
The inclusion $\supseteq$ being clear, we proceed to prove $\subseteq$.
First, $V^{\gen{g}}=V_0^{\gen{g}}\otimes_{k_0}k$, where 
$V_0=A\otimes_{\Z}k_0$ and $k_0$ denotes the prime subfield of $k$. 
If $k_0=\Q$, then clearly
$V_0^{\gen{g}}=A^{\gen{g}}\otimes_{\Z}\Q$. So assume that 
$k_0=\mathbb{F}_p$. Then the $\gen{g}$-cohomology
sequence that is associated with 
$A\stackrel{p}{\mono} A\onto V_0=A/pA$
in conjunction with the fact that $\co1{\gen{g}}{A/pA}$ is trivial
(because $p=\ch k$ does not divide the order of $g$) proves that
$A^{\gen{g}}$ maps onto $V_0^{\gen{g}}$, which finishes the proof.
\end{proof}


\subsection{Singularities}\label{sing}
The \emph{singular locus} of $R$ is defined by
$$
\Sing(R)=\{\p\in\Spec(R)\mid\gldim(R_\p)=\infty\}\ ;
$$
it is a closed subset of $\Spec(R)$ of codimension at least 2
(e.g., \cite[Chapt.~VI]{Ku}).
The complement will be denoted $\Reg(R)$.

\begin{lem} 
Assume $k$ is algebraically closed.
Let $\M$ be a maximal ideal of $S$ and put $\m =\M\cap R$, a 
maximal ideal of $R$. Then $\m\in\Reg(R)$ if and only if
$G^T(\M)$ is a reflection group on $A$.
\end{lem}

\begin{proof}
In view of Lemma (\ref{cotangent}), this is immediate from the
following criterion of Serre \cite{Se} 
(cf.~ also \cite[Exercise 7 on p.~138]{Bou}):
\begin{quote}
$\m\in\Reg(R)$ iff $G^T(\M)$
acts as a pseudoreflection group on $\M/\M^2$.
\end{quote}
\end{proof}

In case $G$ acts without reflections on $A$, the
foregoing leads to a particularly manageable
description of $\Sing(R)$.
For this, we put
$$
I=\bigcap_{1\neq g\in G} I(g)\quad\text{with}\quad
I(g)=(s-s^g \mid s\in S)S\ .
$$
The ideal $I$ is $G$-stable and
semiprime, with
\begin{equation}\label{E:height}
\hgt I = \min_{1\neq g\in G} \rank(1-g)_A
\end{equation}
(see \cite[Lemma 3.2]{L2} and \cite[2.6]{BL}).
So $G$ acts without reflections on $A$ if and only if
$\hgt I\ge 2$.

\begin{cor}
Assume that $G\neq\gen{1}$ acts without reflections on $A$. 
Then, via Lying Over,
$$
\Sing(R)\bij\{\M\in\Spec(S)\mid \M\supseteq I\}/G \ .
$$
This set contains at least two elements.
\end{cor}

\begin{proof}
Recall that Lying Over yields a one-to-one correspondence of
$\Spec(R)$ with $\Spec(S)/G$, the set of $G$-orbits in $\Spec(S)$:
The primes $\p$ of $R$ are exactly the ideals of the form
$\p=\P\cap R$, where $\P$ is a prime of $S$, said to ``lie over"
$\p$, and the full set of primes of $S$ lying over a
particular prime of $R$ forms a $G$-orbit.

Now let $\m$ be a maximal ideal of $R$ and let $\M$ be a maximal
ideal of $S$ lying over $\m$. Then, by the above Lemma,
$\m\in\Reg(G)$ if and only if $G^T(\M)=\gen{1}$. In other words,
since $g\in G^T(\M)$ is equivalent with $\M\subseteq I(g)$,
we have
$$
\m\in\Sing(R)\iff\M\subseteq I\ .
$$

An arbitrary prime $\p$ of $R$ is the intersection of
all maximal ideals $\m\supseteq\p$, and $\p$ belongs to $\Sing(R)$
precisely if all these $\m$'s do. This implies the description of
$\Sing(R)$.

The kernel of the distinguished augmentation of $S=k[A]$ is a
$G$-stable maximal ideal of $S$ containing $I$, and hence it
accounts for a point in $\Sing(R)$. If it was the only point, then
$S/I=k$. But, for any element $g\in G$ of prime order, $\det(1-g)$
is divisible by the same prime, and so $A^{1-g}\neq A$. Therefore,
$S/I(g)\cong k[A/A^{1-g}]\neq k$, and so $S/I\neq k$.
\end{proof}


\subsection{A negative result}\label{negative}
The group $G$ is said to act \emph{fixed point freely} on $A$ if
$A^{\gen{g}}$ is trivial for all $1\neq g\in G$. 
By (\ref{sing})(\ref{E:height}), this
is equivalent with $\hgt I=\rank(A)$, which
in turn just says that $S/I$ is finite $k$-dimensional.
Therefore, as long as $G$ acts without reflections on $A$ and
$k$ is algebraically closed, Corollary (\ref{sing}) implies that
\begin{quote}
\emph{$\Sing(R)$ is finite if and only if $G$ acts fixed point 
freely on $A$.}
\end{quote}
This observation will be used in the proof of the following

\begin{thm} If $G$ acts fixed point freely on $\bar{A}=A/A^G$ and
$\rank(\bar{A})\ge 2$ then $R$ is not a semigroup algebra.
\end{thm}

\begin{proof}
Suppose, by way of contradiction, that $R$ is a 
semigroup algebra. Then, by Proposition 
(\ref{multinv=semigroup}), so is $\bar{R}$. 
Thus we mave assume that $G$ does in fact act fixed 
point freely on $A$. Furthermore, extending scalars 
if necessary, we may assume $k$ to be algebraically 
closed.

As we have remarked above, $\Sing(R)$ is finite. 
On the other hand, $R$ is a semigroup algebra, 
say $R\cong k[M]$. Necessarily, $M$
is an affine normal semigroup with trivial group 
of units, as $A^G$ is trivial. Thus, the action of 
the torus $T$ on $\Max R$, as described in (\ref{torus}), 
has exactly one fixed point.
This action stabilizes $\Sing(R)$.
Since $T$ is connected and $\Sing(R)$ is finite, $T$
must act trivially on $\Sing(R)$.
We conclude that in fact $\#\Sing(R) = 1$, contradicting
Corollary (\ref{sing}). This finishes the proof.
\end{proof}

Since effective lattices for groups of prime order are clearly
fixed point free, we obtain the following 

\begin{cor} 
If $G$ has odd prime order then $R$ is not a semigroup algebra.
\end{cor}


\subsection{An example}\label{example}
If $\rank A=2$ then finite subgroups of $\GLA$ either act
fixed point freely or else they are generated by reflections;
see, e.g., \cite[2.7]{L}. This is of course no longer true in higher
ranks. Here we discuss a specific example which is not directly
covered by the foregoing. Nevertheless, a look at the singularities
very much like the proof of Theorem (\ref{negative})
still yields the desired conclusion. The example is taken from
\cite{KM}, where it was used for different (very interesting) 
purposes.

For a given $n$, let $A=\gen{a_1}\times\ldots\times\gen{a_n}$ be
free abelian of rank $n$, and let 
$G=\diag(\pm 1,\ldots,\pm 1)_{n\times n}\cap\operatorname{SL}(A)$.
So $G$ contains
no reflections but does not act fixed point freely if $n>2$.
We assume $k$ algebraically closed with $\ch k\neq 2$. 

It is not hard to check that the algebra $R$ of multiplicative 
$G$-invariants has the presentation 
$$
R\cong k[x_1,\ldots,x_n,z]/\left(z^2-\prod_{i=1}^n(x_i^2-1)\right)\ .
$$
Using the Jacobian criterion of \cite[p.~173]{Ku}, one checks that
$\Sing R$ is the union of $4\binom{n}{2}$ affine spaces 
$\mathbb{A}_k^{n-2}$, the irreducible components of $\Sing R$.

Alternatively, whenever $G$ has no reflections, such as in our
example, the irreducible
components of $\Sing R$ correspond to the $G$-orbits of the
minimal primes over the ideal $I$ in Corollary (\ref{sing}).
Here, one calculates easily that 
$$
I=\bigcap_{i\neq j\in\{1,\ldots,n\}}P_{i,j}^{\pm,\pm}
\qquad\text{with}\qquad P_{i,j}^{\pm,\pm}=(a_i\pm 1,a_j\pm 1)S \ .
$$
The $P_{i,j}^{\pm,pm}$ are the minimal primes over $I$; they are all
$G$-invariant, and hence they correspond to the irreducible components
of $\Sing R$. Moreover, 
$R/R\cap P_{i,j}^{\pm,\pm}=\left(S/P_{i,j}^{\pm,pm}\right)^G=
k[A_{i,j}]^G$, where $A_{i,j}$ denotes the sublattice of $A$ that is
spanned by all $a$'s except for $a_i$ and $a_j$. Since $G$ acts on
$A_{i,j}$ as the full group $\diag(\pm 1,\ldots,\pm 1)_{n-2\times n-2}$,
it is easy to see that $k[A_{i,j}]^G$ is a polynomial algebra of
dimension $n-2$. Thus the
irreducible components of $\Sing R$ are affine $(n-2)$-spaces.

Now assume $R$ is a semigroup algebra. Then,
by the connectedness argument used in the proof of Theorem 
(\ref{negative}), the torus action considered there stabilizes
all irreducible components of $\Sing R$, and hence also all their
intersections. In the present case, the minimal nonempty intersections
are a collection of $2^n$ points corresponding to the maximal
ideals $(a_1\pm 1,\ldots,a_n\pm 1)S\cap R$ of $R$. Thus all these
are torus fixed points, while there can only be one. This
contradiction shows that, again, $R$ is not a semigroup algebra.





\begin{thebibliography}{99}

\bibitem{Bou}
  N. Bourbaki, \emph{Groupes et alg\`ebres de Lie}, chap. 4--6, 
  Hermann, Paris, 1968.

\bibitem{BL}
  K. A. Brown and M. Lorenz, \emph{Grothendieck groups and higher
  class groups of commutative invariants},
  Contemp. Math. \textbf{184} (1995), 59--74.

\bibitem{BH} 
  W. Bruns and J. Herzog, \emph{Cohen-Macaulay Rings},
  Cambridge University Press, Cambridge, 1993.

\bibitem{Fa}
  D. R. Farkas, \emph{Multiplicative invariants},
  L'Enseignement Math. \textbf{30} (1984), 141--157.

\bibitem{Fa2}
  \bysame, \emph{Reflection groups and multiplicative invariants},
  Rocky Mt. J. \textbf{16} (1986), 215--222.

\bibitem{Fa3}
  \bysame, \emph{Toward multiplicative invariant theory},
  in:``Group Actions on Rings" (S. Montgomery, ed.), Amer. Math. Soc.,
  Providence, RI, Contemporary Math. \textbf{43} (1985), 69--80. 

\bibitem{F}
  R. M. Fossum, \emph{The Divisor Class Group of a Krull
  Domain}, Springer-Verlag, Berlin-Heidelberg, 1973.

\bibitem{Fu}
  W. Fulton, \emph{Introduction to Toric Varieties},
  Annals of Mathematics Studies, No. 131, 
  Princeton University Press, Princeton, 1993

\bibitem{G}
  R. Gilmer, \emph{Commutative Semigroup Rings},
  Chicago Lectures in Mathematics, 
  The University of Chicago Press,
  Chicago and London, 1984.

\bibitem{Gu}
  J. Gubeladze, \emph{Anderson's conjecture and the maximal
  monoid class over which projective modules are free},
  Math. USSR Sbornik \textbf{63} (1989), 165--180.
 
\bibitem{Gu2}
  \bysame, \emph{The elementary action on unimodular
  rows over a monoid ring},
  J. Algebra \textbf{148} (1992), 135--161.


\bibitem{H}
  J. E. Humphreys, \emph{Reflection Groups and Coxeter Groups},
  Cambridge University Press, Cambridge, 1990.

\bibitem{KM}
  P. H. Kropholler and B. Moselle, \emph{A family of 
  crystallographic groups with $2$-torsion if $K_0$ of the 
  rational group algebra},
  Proc. Edinburgh Math. Soc. \textbf{34} (1991), 325--331.

\bibitem{Ku}
  E. Kunz, \emph{Introduction to Commutative Algebra and
  Algebraic Geometry},
  Birkh\"auser, Boston, 1985.

\bibitem{Lem}
  N. Lemire, \emph{Multiplicative Invariants},
  Ph.D. thesis, University of Alberta, Edmonton, 1998.

\bibitem{L}
  M. Lorenz, \emph{Class groups of multiplicative invariants},
  J. Algebra \textbf{177} (1995), 242--254.

\bibitem{L1}
  \bysame, \emph{Regularity of multiplicative invariants},
  Comm. Algebra \textbf{24} (1996), 1051--1055.

\bibitem{L2}
  \bysame, \emph{Picard groups of multiplicative invariants},
  Comment. Math. Helv. \textbf{72} (1997), 389--399.
  
\bibitem{Sch}
  M. Sch\"onert et al., {\sf GAP} -- \emph{Groups, Algorithms 
  and Programming}, Lehrstuhl D f{\"u}r Mathematik,
  Rheinisch Westf{\"a}lische Technische Hoch\-schule,
  Aachen, Germany, fifth edition, 1995; 
  available via anonymous ftp from \texttt{math.rwth-aachen.de}.

\bibitem{Se}
  J.-P. Serre, \emph{Groupes finis d'automorphismes d'anneaux 
  locaux r\'eguliers}, in: Colloque d'al\-g\`eb\-re, No. 8,
  \'Ecole Norm. Sup. des Jeunes Filles, Paris, 1967.

\bibitem{St}
  B. Sturmfels, \emph{Gr\"obner Bases and Convex Polytopes},
  University Lecture Series, Vol. 8, Amer. Math. Soc., 
  Providence, 1996.

\end{thebibliography}
\end{document}